\documentclass[12pt,draft]{article}
\usepackage{amsmath,amsthm,amsfonts,amssymb,enumerate}
\date{}
\newtheorem{Theorem}{Theorem}[section]

\newtheorem{Proposition}[Theorem]{Proposition}
\newtheorem{Lemma}[Theorem]{Lemma}
\newtheorem{Corollary}[Theorem]{Corollary}
\newtheorem{Remark}[Theorem]{Remark}

\numberwithin{equation}{section}

\def\R{\mathbb R}
\def\N{\mathbb N}

\def\E{\mathbb E}
\def\P{\mathbb P}
\def\M{\mathbb M}

\begin{document}

\begin{center}
{\bf\Large  Weak solutions to  the  stochastic porous media equation via Kolmogorov equations: the degenerate case}
\vspace{5mm}

Viorel Barbu (University of Iasi, 6600 Iasi, Romania).\vspace{5mm}

Vladimir I. Bogachev (Department of Mechanics and Mathematics,
         Moscow State University, 119899 Moscow, Russia).\vspace{5mm}

Giuseppe Da Prato (Scuola Normale Superiore, Piazza dei Cavalieri 7, 56126 Pisa,  Italy)
\vspace{5mm}

Michael R\"ockner (Fakult\"at f\"ur Mathematik, Universit\"at   Bielefeld,
Postfach 100131, D--33501 Bielefeld,
Germany)
\end{center}\vspace{5mm}

{\bf Abstract}. A stochastic version of the porous medium equation
with coloured noise is studied. The corresponding Kolmogorov equation
is solved in the space $L^2(H,\nu)$ where $\nu$ is an
infinitesimally excessive
measure. Then a weak solution is constructed.

\vskip .1cm

\noindent {\bf 2000 Mathematics Subject Classification AMS}:
76S05, 35J25, 37L40.

\section{Introduction}

The porous medium equation
\begin{equation}
\label{e1.1}
\frac{\partial X}{\partial t}=\Delta(\Psi(X)),\quad m\in \N,
\end{equation}
on a bounded open set $D\subset \R^d$ with Dirichlet boundary
conditions for the Laplacian $\Delta$ and with $\Psi$ in a large class
of functions has been studied extensively (see e.g. \cite{A86},
\cite[Section 4.3]{Ba93}). Recently, there has been also several
papers on the stochastic version of \eqref{e1.1}, i.e.

\begin{equation}
\label{e1.2}
dX(t)= \Delta(\Psi(X(t))dt+\sqrt{C}\;dW(t),\quad t\ge 0,
\end{equation}
(cf. \cite{DPR03a}, \cite{DPR03b}, \cite{BDPR03} and \cite{BDP02}).

In this paper we continue the study of the stochastic partial
differential equation (SPDE) \eqref{e1.2}. Before we describe our new
results precisely, let us fix some notation and our exact conditions.

The appropriate state space on which we consider \eqref{e1.2} is
$H:=H^{-1}(D)$, i.e., the dual of the Sobolev space
 $H^1_0:=H^1_0(D)$, with inner product $\langle\cdot,\cdot \rangle_H.$
Below we shall use the standard $L^2(D)$-dualization $\langle
\cdot,\cdot \rangle_H$ between $H^1_0(D)$ and $H=H^{-1}(D)$ induced by
the embeddings

$$
H^1_0(D)\subset L^2(D)'=L^2(D)\subset H^{-1}(D)=H
$$
  without further notice. Then for $x\in H$ one has
$$
|x|^2_H=\int_D((-\Delta)^{-1}x)(\xi)\;x(\xi)d\xi.
$$

Let
$(W_t)_{t\ge 0}$ be a cylindrical Brownian motion in $H$ and let
$C$ be a positive definite bounded operator on $H$  of
trace class. To be more concrete below we assume:

\begin{itemize}
\item[$(H1)$]
{\it There exist numbers $\lambda_k\in[0,\infty)$,
$k\in \N$ such that for the
eigenbasis $\{e_k|\;k\in \N\}$ of $\Delta$
in $H$ {\rm(}with Dirichlet boundary
conditions{\rm)} we have
$Ce_k= \lambda_ke_k$ for all $k\in \N$.}

\item[$(H2)$]
     {\it
For $\alpha_k:=\sup_{\xi\in D}|e_k(\xi)|^2,\; k\in \N$, we have
$K:=\sum_{k=1}^\infty\alpha_k\lambda_k<+\infty.$}

\item[$(H3)$]
{\it
There exist $\Psi\in C^1(\R), r\in (1,\infty), \kappa_0, \kappa_1,
 C_1>0$ such that}
$$
\kappa_0|s|^{r-1}\le \Psi'(s)\le \kappa_1|s|^{r-1}+C_1
\quad \hbox{for all $s\in \R$ {\rm(}cf. {\rm\cite{BDPR03})}.}
$$
\end{itemize}

Our general aim in studying SPDE \eqref{e1.2} is to construct a strong
Markov weak solution for \eqref{e1.2}, i.e., a solution in the sense
of the corresponding martingale problem (see \cite{SV79} for the
finite dimensional case), at least for a large set $\overline{H}$ of
starting points in $H$ which is invariant for the process, i.e., with
probability one $X_t\in \overline{H}$ for all $t\ge 0.$ We follow the
strategy first presented in \cite{R99} (and already carried out in
cases with bounded $C^{-1}$ in \cite{DPR02}). That is, first we want
to construct a solution to the corresponding Kolmogorov equations in
$L^2(H,\mu)$ for suitably chosen reference measures $\mu$ (see below),
and then a strong Markov process with continuous sample paths having
transition probabilities given by that solution to the Kolmogorov
equations. As in \cite{DPR02} we also aim to prove that this process
is for $\mu$-a.e. starting point $x\in \overline{H}$ a unique (in
distribution) continuous Markov process whose transition semigroup
consists of continuous operators on $L^2(H,\mu)$, which is e.g. the
case if $\mu$ is one of its excessive measures.

Before we summarize the specific new results of this paper in relation
to those obtained in \cite{DPR03a},\cite{DPR03b},
\cite{BDPR03}, let us describe this programme more precisely.

\vspace{2mm}

Applying It\^o's formula (at a heuristic level) to \eqref{e1.2} one
finds what the corresponding Kolmogorov
operator, let us call it $N_0$, should be, namely
\begin{equation}
\label{e1.3}
N_0\varphi(x)=\frac12\;
\sum_{k=1}^\infty\lambda_kD^2\varphi(e_k,e_k)
+D\varphi(x)(\Delta(\Psi(x))),\quad x\in H,
\end{equation}
where $D\varphi$, $D^2\varphi$ denote the first and second Fr\'echet
derivatives of $\varphi:H\to\R$. So, we take $\varphi\in C^2_b(H)$.

In order to make sense of \eqref{e1.3} one needs that
$\Delta(\Psi(x))\in H$ at least for ``relevant'' $x\in H$. Here one
clearly sees the difficulties since $\Psi(x)$ is, of course, not
defined for any Schwartz distribution in $H=H^{-1}$, not to mention
that it will not be in $H^1_0(D).$ So, a way out of this is to think
about ``relevant'' $x\in H.$ Our approach to this is first to look for
an invariant measure for the solution to equation \eqref{e1.2} which
can now be defined ``infinitesimally''(cf. \cite{BR01}) without having
a solution to \eqref{e1.2} as a solution to the equation
\begin{equation}
\label{e1.4}
N_0^*\mu=0
\end{equation}
with the property that $\mu$ is supported by those $x\in H$ for which
$\Psi(x)$ makes sense and $\Delta(\Psi(x))\in H$.
Equation \eqref{e1.4} is a short form for
\begin{equation}
\label{e1.5}
N_0\varphi\in L^1(H,\mu)\quad \hbox{and}\quad
\int_HN_0\varphi d\mu=0\quad \hbox{for all $\varphi\in C_b^2(H)$.}
\end{equation}
Any invariant measure for any solution of \eqref{e1.2} in the classical sense will satisfy \eqref{e1.4}.
Then we can
analyze $N_0$, with domain $C_b^2(H)$ in $L^2(H,\mu)$, i.e.,
solve the Kolmogorov equation
\begin{equation}
\label{e1.6}
\frac{dv}{dt}=\overline{N_0}v,\quad v(0,\cdot)=f,
\end{equation}
for the closure $\overline{N_0}$ of $N_0$ on $L^2(H,\mu)$ and initial
condition $f\in L^2(H,\mu)$. This means, we have to prove that
$\overline{N_0}$ generates a $C_0$-semigroup
$T_t=e^{t\overline{N_0}}$ on $L^2(H,\mu)$, i.e., that $(N_0,C_b^2(H))$
is essentially $m$-dissipative on $L^2(H,\mu)$.

Subsequently, we have to show that $(T_t)_{t\ge 0}$ is given by a
semigroup of probability kernels $(p_t)_{t\ge 0}$ (i.e., $p_tf$ is a
$\mu$-version of
 $T_tf\in L^2(H,\mu)$ for any $t\ge 0$ and any
bounded measurable function $f\colon\, H\to \R$)
and such that there exists a strong Markov process with
continuous sample paths in $H$ whose transition function is
$(p_t)_{t\ge 0}.$ Then, by definition, this Markov process will solve
the martingale problem corresponding to \eqref{e1.2}.

The existence of solutions to \eqref{e1.4} (even for more general SPDE
than \eqref{e1.2}) was proved in \cite{BDPR03} (the method was based
essentially on finite dimensional approximations), generalizing earlier
results from \cite{DPR03a}. We shall restate the precise
theorem in \S 2 below.

In \cite{DPR03a} in the special case when
\begin{equation}
\label{e1.7}
\Psi(s):=\alpha s+s^m,\quad s\in \R,
\end{equation}
for $m\in \N$, $m$ odd, and $\alpha>0$, the remaining part of the above programme was carried out. The specially
interesting ``degenerate'' case $\alpha= 0$ in \eqref{e1.7} was, however, not covered.\bigskip

In this paper we shall improve these results in an essential way.
First, we shall construct a solution to the Kolmogorov equation
\eqref{e1.6} for $\Psi$ as in (H3), hence including the case
$\alpha=0$ in \eqref{e1.7}. More precisely, we identify a whole class
$\mathcal M$ of reference measures,
called infinitesimally excessive measures,
which includes all measures solving
\eqref{e1.4} so that for all $\nu\in \mathcal M$ we can construct a
solution to the Kolmogorov equation \eqref{e1.6} in $L^2(H,\nu)$ for
$\Psi$ as in (H3), hence including the degenerate case $\alpha=0$, in
\eqref{e1.7}. The main tool employed here is the Yosida approximation
for the nonlinear maximal dissipative mapping $\Delta(\Psi)$, as a map
in $H^{-1}$ with domain $H^1_0$. In particular, we thus clarify that in
case the nonlinearity of SPDE \eqref{e1.2} is maximal dissipative, the
issue of proving the existence of infinitesimally invariant measures $\mu$
for $N_0$ and the issue of essentially maximal dissipativity of the
operator $(N_0,C^2_b(H))$ on $L^2(H,\nu)$ can be separated completely.
That is, the latter does not depend in particular on how one
constructs a solution to \eqref{e1.4} and which solution is chosen as
a reference measure.

Second, we shall construct the said Markov process which weakly solves
SPDE \eqref{e1.2} for general $\Psi$ as in (H3); i.e., without any
nondegeneracy assumptions. Furthermore, we prove that for $d=1$ and
specifically chosen $C$ (cf. condition (H.4) in \S 5) the Markov
process is strong Feller.

The organization of this paper is as follows. In \S 2 we summarize all
relevant results about infinitesimal invariant measures $\mu$ for
$N_0$ from \cite{BDPR03} and \cite{DPR03a}. Then we define the
mentioned class $\mathcal M$ of references measures $\nu$ and prove
that for some $\lambda>0$, $(\lambda-N_0,C^2_b(H))$ is dissipative on
$L^2(H,\nu)$, hence $(N_0,C^2_b(H))$ is closable on $L^2(H,\nu)$.

\S 3 is devoted to the Yosida approximations. In \S 4 we prove that
for all $\nu\in \mathcal M$ the closure of $(N_0,C^2_b(H))$ on
$L^2(H,\nu)$ generates a $C_0$-semigroup on $L^2(H,\nu)$ solving
\eqref{e1.6}. \S 5 is devoted to the existence and uniqueness of a Markov
process solving SPDE \eqref{e1.2} in the sense of a martingale
problem, and, in case $d=1$, to its strong Feller property on $\mathop{supp}\nu$.
In \S 6 under weak additional conditions we prove that
if $\nu$ is the solution of \eqref{e1.4} constructed in
\cite{BDPR03}, then $\mathop{supp}\nu = H$, i.e.\ $\nu$ charges any 
non-empty open set of $H$.

Finally, we would like to mention that we think that it should be also
possible to prove the existence and uniqueness of a strong solution for
(1.2). A corresponding paper of the last named author jointly with B.
Rozovskii is in preparation.

\section{Infinitesimal invariance  and a large class of references measures}

We first note that
$N_0\varphi(x)$ is well defined for $\varphi\in C^2_b(H)$ if $x$ belongs to the set
\begin{equation}
\label{e2.1}
H_\Psi:=\{x\in L^2(D):\; \Psi(x)\in H^1_0  \}.
\end{equation}
We also define for $r>1$
$$
H_{0,r}^1:=\{x\in L^2(D):\; |x|^r\mbox{\rm sign}\;x \in H^1_0  \}.
$$
Now we recall the following result from \cite{BDPR03} (see Theorem 1.1 and Corollary 1.1 ibid).
\begin{Theorem}
\label{t2.1}
Assume that $(H1)$--$(H3)$ hold. Then there exists a probability measure $\mu$ on $H$ which is
infinitesimally invariant for $N_0$ in the sense of \eqref{e1.5}. Furthermore,
\begin{equation}
\label{e2.2}
\int_H\int_D|\nabla(\Psi(x))(\xi)|^2d\xi \mu(d x)<+\infty
\end{equation}
and
\begin{equation}
\label{e2.3}
\int_H\int_D|\nabla(|x|^{\frac{r+1}{2}}\;\mbox{\rm sign}\;x)(\xi)|^2d\xi \mu(d x)<+\infty.
\end{equation}
In particular, $\mu({H_\Psi\cap H^1_{0,\frac{r+1}{2}}})=1.$
\end{Theorem}
\begin{Remark}
\label{r2.2}
{\em It was also shown in \cite[Lemma 1]{BDPR03} that $H_\Psi\subset H^1_{0,r}.$
So, \eqref{e2.2} implies that
$$
\int_H\int_D|\nabla(|x|^{r}\;\mbox{\rm sign}\;x)(\xi)|^2d\xi \mu(d x)<+\infty.
$$
Therefore, by Poincar\'e's inequality, $H_\Psi\subset L^{2r}(D)$ and
$$
\int_H\int_D |x|^{2r}(\xi)d\xi \mu(d x)<+\infty.
$$}
\end{Remark}

By Theorem \ref{t2.1}, $N_0\varphi$ is $\mu$-a.e. defined for all
$\varphi\in C^2_b(H)$. All subsequent results in this paper are valid
for the larger class of measures $\mathcal M$ on $H$ which contains all
infinitesimally invariant measures for $N_0$ and
consists of all probability measures $\nu$ on $H$ which satisfy
\eqref{e2.2} and for which there exists $\lambda_\nu\in (0,\infty)$
such that

\begin{equation}
\label{e2.4}
\int_HN_0\varphi d\nu\le \lambda_\nu\int_H\varphi d\nu
\quad\hbox{for all $\varphi\in C^2_b(H)$ with $\varphi\ge 0$
\ $\nu$-a.e.}
\end{equation}
The elements in $\mathcal M$ are called 
\emph{infinitesimally excessive measures}.
\begin{Lemma}
\label{l2.3}
Let $\nu\in \mathcal M$ and $\varphi\in C^2_b(H)$ be
such that $\varphi=0$\ $\nu$-a.e.
Then $N_0\varphi=0$\ $\nu$-a.e.
\end{Lemma}
\begin{proof}
 The proof is analogous to that of
Lemma 3.1 in \cite{BDPR03} (see also \cite[Proposition 4.1]{RS03b}).
\end{proof}

We would like to emphasize that so far we have not been able to show
that $\mu(U)>0$ (for $\mu$ as in Theorem \ref{t2.1}) for any open non
empty set $U\subset H$. So, Lemma \ref{l2.3} is crucial to consider
$N_0$ as an operator on $L^2(H,\mu)$ with domain equal to the
$\mu$-classes determined by $C^2_b(H)$, again denoted by $C^2_b(H)$.
The same holds for any $\nu\in \mathcal M$.

\begin{Proposition}
\label{p2.4}
Assume that $(H1)$--$(H3)$ hold and let $\nu\in \mathcal M$. Then:
\begin{enumerate}
\item[{\rm(i)}] $N_0\varphi\in L^2(H,\nu)$ for $\varphi\in C^2_b(H)$.

\item[{\rm(ii)}] $(\frac12\lambda_{\nu}-N_0,C^2_b(H))$ is dissipative
on $L^2(H,\nu)$,
i.e.,
$$
\|\lambda(\lambda+
\frac12\lambda_{\nu}-N_0\varphi\|_{L^2(H,\nu)}\ge \|\varphi\|_{L^2(H,\nu)}
\ \hbox{for all $\varphi\in C^2_b(H)$.}
$$
In
particular,
$(N_0,C^2_b(H))$ is closable on $L^2(H,\nu)$, its closure being denoted by $(N_2,D(N_2))$.
\end{enumerate}
\end{Proposition}
\begin{proof}
(i) We note that
$$
\int_D|\nabla \Psi(x)|^2(\xi)d(\xi)=|\Delta\Psi(x)|^2_H.
$$
Hence the assertion follows by \eqref{e2.2}.

(ii) This follows from \cite[Appendix B, Lemma  1.8]{E99}.
\end{proof}

\section{Yosida approximations}

For completeness we recall the definition and basic properties of
the Yosida approximation of an
$m$-dissipative
map $F\colon D(F)\subset H\to H$. The latter means that
\begin{equation}
\label{e3.1}
\langle F(x)-F(y),x-y  \rangle_H\le 0\quad\mbox{\rm for all}\;x,y\in D(F)
\end{equation}
and
\begin{equation}
\label{e3.2}
(\lambda I-F)(D(F))=H \quad\mbox{\rm for all}\;\lambda>0,
\end{equation}
where $I$ denotes the identity operator on $H$.
For $\varepsilon>0$ let
\begin{equation}
\label{e3.3}
J_\varepsilon:=(I-\varepsilon F)^{-1}.
\end{equation}
Note that by \eqref{e3.1} $I-\varepsilon F:D(F)\to H$ is one--to--one.
Then $J_\varepsilon$  is Lipschitz
continuous with constant $1$, hence so is
\begin{equation}
\label{e3.4}
F_\varepsilon:=\frac1\varepsilon(J_\varepsilon-I)
\end{equation}
with constant $\frac{1}{\varepsilon}$.
 $F_\varepsilon$ is called {\em Yosida approximation}
of $F$. It has the following properties (cf. e.g. \cite{Ba93},   \cite{DP76} or \cite{Sho97}):
\begin{equation}
\label{e3.5}
\lim_{\varepsilon\to 0}F_\varepsilon(x)=F(x),\quad x\in D(F),
\end{equation}
\begin{equation}
\label{e3.6}
|F_\varepsilon(x)|_H\le |F(x)|_H,,\quad x\in D(F),\;\varepsilon>0,
\end{equation}
\begin{equation}
\label{e3.7}
|F_\varepsilon(x)|_H\uparrow 1_{D(F)}(x) |F(x)| +\infty\cdot 1_{H\setminus D(F)}(x),\quad\mbox{\rm
as}\;\varepsilon\downarrow 0,\;x\in H,
\end{equation}
\begin{equation}
\label{e3.8}
\langle F_\varepsilon(x),F(x)\rangle_H\le -|F_\varepsilon(x)|^2_H,\quad x\in D(F).
\end{equation}
The following is well known, see e.g. \cite[Chapter 2, Proposition 2.12]{Ba93}
and for the original proof \cite{Brezis71}.

\begin{Proposition}
\label{p3.1}
Assume $(H3)$ holds. Then $F:=\Delta\Psi$ with domain $D(F)\colon =H_\Psi$ is
$m$-dissipative on $H$.
\end{Proposition}

\section{Essential maximal dissipativity of $N_0$ on $L^2(H,\nu)$}

Below, $F_\varepsilon,\;\varepsilon>0$, shall always denote the Yosida approximation to
$(\Delta\Psi, H_\Psi)$. We need a further regularization and, therefore, define for $\beta>0$
\begin{equation}
\label{e4.1}
 F_{\varepsilon,\beta}(x):=\int_H e^{\beta B}F_\varepsilon(e^{\beta B}x+y)N_{\frac12\;B^{-1}(e^{2\beta B}x-I)},
\quad x\in H,
\end{equation}
where $B\colon D(B)\subset H\to H$ is a self-adjoint negative definite operator such that $B^{-1}$ is of trace
class. Then obviously $F_{\varepsilon,\beta}$ is dissipative of class $C^\infty$, and has bounded derivatives of all
orders. Furthermore,
\begin{equation}
\label{e4.2}
 \lim_{\beta\to 0}F_{\varepsilon,\beta}(x)=F_{\varepsilon}(x),
\quad x\in H,
\end{equation}
(see \cite[Theorem 9.19]{DPZ1}) and, since $F_{\varepsilon}$ is Lipschitz, there exists $c_\varepsilon\in
(0,\infty)$ such that
\begin{equation}
\label{e4.3}
|F_{\varepsilon,\beta}(x)|\le c_\varepsilon(1+|x|_H),\quad x\in H,
\;
\quad x\in H,\beta>0.
\end{equation}
Now consider the approximating stochastic equation
\begin{equation}
\label{e4.4}
dX(t)=F_{\varepsilon,\beta}(X(t))dt+\sqrt{C}\;dW(t)
\end{equation}
It is well known  (see \cite{DPZ1}) that for any initial condition $x\in H$ equation \eqref{e4.4}
has a unique solution $X_{\varepsilon,\beta}(\cdot,x)$ and that for $\lambda>0$ and $f\in C^2_b(H)$
\begin{equation}
\label{e4.5}
\varphi_{\varepsilon,\beta}(x)=\int_0^\infty e^{-\lambda t}\E[f(X_{\varepsilon,\beta}(t,x))]dt,\quad x\in H,
\end{equation}
is in $C^2_b(H)$ and solves the equation
\begin{equation}
\label{e4.6}
f(x)=\lambda \varphi_{\varepsilon,\beta}(x)-\frac12\;\sum_{k=1}^\infty\lambda _kD^2\varphi_{\varepsilon,\beta}(x)
(e_k,e_k)+ D\varphi_{\varepsilon,\beta}(x)(F_{\varepsilon,\beta}(x)),
\end{equation}
(see \cite[Chapter 5.4]{DPZ2}).
We have moreover
for all $h\in H,$
\begin{equation}
\label{e4.7}
D\varphi_{\varepsilon,\beta}(x)(h) =\int_{0}^{+\infty}e^{-\lambda
 t}\E[ Df(X_{\varepsilon,\beta}(t,x))(D_xX_{\varepsilon,\beta}(t,x)h)]dt.
\end{equation}
For any  $h\in H$ we  set $\eta_{\varepsilon,\beta
}^h:=D_xX_{\varepsilon,\beta}(t,x).$ Then we have
\begin{equation}
\label{e4.8}
\left\{\begin{array}{l}
\frac{d}{dt}\;\eta_{\varepsilon,\beta }^h(t,x)
= DF_{\varepsilon,\beta }(X_{\varepsilon,\beta
}(t,x))\eta_{\varepsilon,\beta }^h(t,x) \\\\
\eta_{\varepsilon,\beta }^h(0,x)=h.
\end{array}\right.
\end{equation}
Multiplying both sides of   equation \eqref{e4.8} by $\eta_{\varepsilon,\beta
}^h(t,x),$ integrating with respect to $t$ and taking  the
dissipativity of  $DF_{\varepsilon,\beta }$ into account,  we find

\begin{equation}
\label{e4.9}
|\eta_{\varepsilon,\beta }^h(t,x)|^2\le |h|^2.
\end{equation}
Consequently by \eqref{e4.7} it follows that
\begin{equation}
\label{e4.10}
 |D\varphi_{\varepsilon,\beta}(x)|_{H^1_0} \le \frac{1}{\lambda }\;\|Df\|_{0},\;x\in H,
\end{equation}
where $\|\cdot\|_{0}$ denotes the sup norm.

Now we can prove the following result.

\begin{Theorem}
\label{t4.1}
 Assume that $(H1)$--$(H3)$ hold and let $\nu\in \mathcal M$. Then $(N_0,C^2_b(H))$  is essentially
$m$-dissipative
on $L^2(H,\nu),$ i.e., its closure $(N_2,D(N_2))$ is maximal dissipative on $L^2(H,\nu)$.
\end{Theorem}
\begin{proof}
Let $f\in C^2_b(H)$ and let $\varphi_{\varepsilon,\beta}$ be the
solution to equation \eqref{e4.6}. Then
$\varphi_{\varepsilon,\beta}\in C^2_b(H)$ and we  have
\begin{equation}
\label{e4.11}
\lambda \varphi_{\varepsilon,\beta}-N_0\varphi_{\varepsilon,\beta}
=f+ D\varphi_{\varepsilon,\beta}(F_{\varepsilon,\beta}-\Delta\Psi).
\end{equation}
We claim that
$$\lim_{\varepsilon \to 0}\lim_{\beta \to 0} D\varphi_{\varepsilon,\beta}(F_{\varepsilon,\beta}-\Delta\Psi)=0\;\;\mbox{\rm
in}\;L^2(H,\nu).
$$
In fact by \eqref{e4.10} it follows that
\begin{equation}\label{e4.12}
I_{\varepsilon,\beta }:=
\int_{H}^{}|D\varphi_{\varepsilon,\beta}
(F_{\varepsilon,\beta}-\Delta\Psi)|_{H^1_0}^2d\nu \le
\frac{1}{\lambda^{2} }\;\|Df\|^{2}_{0}
\int_{H}^{}|F_{\varepsilon,\beta}-\Delta\Psi|_H^2d\nu .
\end{equation}
Letting $\beta\to 0$ we conclude by \eqref{e4.3} that
$$
\limsup_{\beta \to 0}I_{\varepsilon ,\beta }\le\frac{1}{\lambda^2 }\;\|Df\|^2_{0}
\int_{H}^{}|F_{\varepsilon}-\Delta\Psi|_H^2d\nu.
$$
Since $\nu$ verifies  \eqref{e2.2} by assumption, the claim now follows, in view of the dominated convergence
theorem, from
\eqref{e3.6}--\eqref{e3.7} with $F:=\Delta\Psi$.

Hence we  have proved that
$$
\lim_{\varepsilon\to 0}\lim_{\beta \to 0}(\lambda -N_0)\varphi_{\varepsilon,\beta}=f\;\;\mbox{\rm in}\;L^2(H,\nu).
$$
Therefore the closure of the range of $\lambda -N_0$ includes $C^2_b(H)$
which is dense in $L^2(H,\nu).$ By the
Lumer--Phillips theorem it follows that $N_2$ is maximal--dissipative as
required.
\end{proof}

As a consequence of the proof of Theorem \ref{t4.1} we have:

\begin{Corollary}
\label{c4.2}
 Assume that $(H1)$--$(H3)$ hold and let $\nu\in \mathcal M$. Define a
$C_0$-semi\-group
$$
P_t:=e^{tN_2},\quad t\ge 0,
$$
 on $L^2(H,\nu)$ {\rm(}which exists by Theorem {\rm\ref{t4.1})}. Then

\begin{enumerate}
\item[{\rm(i)}]
$v(t,\cdot):=P_tf$, $t>0$,
solves \eqref{e1.6} for the initial datum $f\in D(N_2)$.

\item[{\rm(ii)}] $(P_t)_{t\ge 0}$ is Markovian, i.e., $P_t1=1$ and $P_tf\ge 0$  for all nonnegative $f\in L^2(H,\nu)$
and all $t\ge 0$.

\item[{\rm(iii)}] Let $f\in L^2(H,\nu)$, $f$ nonnegative, and $t>0$. Then
\begin{equation}
\label{e4.13}
\int_HP_tfd\nu\le e^{\lambda_\nu t}\int_Hfd\nu.
\end{equation}
\end{enumerate}
\end{Corollary}
\begin{proof}
(i) The assertion follows by the definition of $P_t,\;t\ge 0$.

(ii) By  \cite[Appendix B, Lemma 1.9]{E99}
$P_t$ is positivity preserving.
Since $1\in C^2_b(H)$ and $N_01=0,$ it follows that $P_t1=1$
for all $t\ge 0$.

(iii) We first note that since $C^2_b(H)$ is dense in $D(N_2)$ with respect to the graph norm given by $N_2$, it
follows by Theorem \ref{t4.1} and \eqref{e2.4} that
\begin{equation}
\label{e4.14}
\int_HN_2fd\nu\le  \lambda_\nu  \int_Hfd\nu,\quad
\hbox{for all $f\in D(N_2)$ with $f\ge 0$\ $\nu$-a.e.}
\end{equation}
So, if $f\in C^2_b(H)$ ($\subset D(N_2)$), $f\ge 0$,
then $P_tf\in D(N_2)$ and $P_tf\ge 0$\ $\nu$-a.e.
Hence \eqref{e4.14} and assertion (i) imply that
$$
\frac{d}{dt}\;\int_HP_tfd\nu =\int_HN_2P_tfd\nu
\le \lambda_\nu\int_HP_tfd\nu.
$$
So, by Gronwall's lemma \eqref{e4.13} follows for $f\in C^2_b(H)$,
$f\ge 0$. But since any nonnegative $f\in L^2(H,\nu)$
can be approximated by
nonnegative functions in $C^2_b(H)$ in $L^2(H,\nu)$, assertion (iii)
follows.
\end{proof}

\section{Existence of a weak solution of SPDE \eqref{e1.2}}

This section generalizes all results of \S 4 in \cite{DPR03b} in an
essential way. However, parts of it are very similar. We,
nevertheless, include a complete presentation below for the
reader's convenience.

\begin{Theorem}[Existence]
\label{t5.1}
 Assume that $(H1)$--$(H3)$ hold and, in addition, that $r\ge 2$.
Let $\nu\in \mathcal M$. Then
\begin{enumerate}
\item[{\rm(i)}] There exists a conservative strong Markov process
$$
\M=\Bigl(\Omega,\mathcal F, (\mathcal F_t)_{t\ge 0},
(X_t)_{t\ge 0},(\P_x)_{x\in H}\Bigr)
$$
on $H$ with continuous sample paths such that for its transition semigroup
$(p_t)_{t\ge 0}$ defined by
$$
p_tf(x):=\int_Hf(X_t)d\P_x,\quad t\ge 0, x\in H,
$$
where $f\colon\, H\to \R$ is
bounded $\mathcal B(H)$-measurable, we have that $p_tf$ is a
$\nu$-version of $e^{tN_2}f$, $t>0$.
Furthermore, if $f\ge 0$, one has
$$
\int_Hp_tfd\nu\le e^{\lambda_\nu t}\int_Hfd\nu\quad
\hbox{for all $t\ge 0$,}
$$
i.e., $\nu$ is an excessive measure for $\M$.

\item[{\rm(ii)}]
There exists $\overline{H}\in \mathcal B(H)$ such that $\nu(\overline{H})=1$,
for all $x\in \overline{H}$  one has
$$
\P_x[X_t\in \overline{H}\quad\forall\;t\ge 0]=1,
$$
and for any probability measures $\rho$ on $(H,\mathcal B(H))$
with $\rho(\overline{H})=1$, the process
$$
\varphi(X_t)-\int_0^tN_0\varphi(X_s)ds,\quad t\ge 0,
$$
is an $(\mathcal F_t)$-martingale under $\P_\rho:=\int_{\overline{H}}\P_x\rho(dx)$
for all $\varphi\in C^2_b(H)$ and one has
$\P_\rho\circ X_0^{-1} =\rho$.
\end{enumerate}
\end{Theorem}

\begin{Theorem}[Uniqueness]
\label{t5.2}
 Assume that $(H1)$--$(H3)$ hold and,   in addition,
that $r\ge 2.$ Let $\nu\in \mathcal M$. Suppose that
$$
\M'=\Bigl(\Omega',\mathcal F', (\mathcal F'_t)_{t\ge 0},(X'_t)_{t\ge 0},
(\P'_x)_{x\in H}\Bigr)
$$
is a continuous Markov process on $H$ whose transition semigroup
$(p'_t)_{t\ge 0}$ consists of continuous operators on $L^2(H,\mu)$
with locally {\rm(}in $t${\rm)} uniformly bounded operator norm
$($which is e.g. the case if $\nu$ is also an excessive measure for
$\M'$$)$. If $\M'$ satisfies assertion {\rm(ii)} of Theorem
{\rm\ref{t5.1}} for $\rho:=\nu$, then for $\nu$-a.e. $x\in H$, one has
$p'_t(x,dy)=p_t(x,dy)$ for all $t\ge 0$ {\rm(}where $p_t$ is as in
Theorem {\rm\ref{t5.1}-(i))}, i.e., $\M'$ has the same finite
dimensional distributions as $\M$ for $\nu$-a.e. starting point.
\end{Theorem}

We shall only prove Theorem \ref{t5.1}-(i). The remaining parts are
proved in exactly the same way as Theorem 7.4-(ii), Proposition 8.2
and Theorem 8.3 in \cite{DPR02} with the only exception that because
we do not know whether $(p_t)_{t\ge 0}$ is Feller, all statements can
only be proved $\nu$-a.e. So we do not want to repeat them here.

Our proof of Theorem \ref{t5.1}-(i) is based on the theory of
generalized Dirichlet forms developed in \cite{st99}. Indeed, by
Corollary \ref{c4.2}, $(N_2,D(N_2))$ is a Dirichlet operator in the
sense of \cite{MR92}, \cite{st99}. Hence
by \cite[Proposition I.4.6]{st99}

$$
\mathcal E(u,v):=\left\{\begin{array}{l}
(u,v)_{L^2(H,\nu)}-(N_2u,v)_{L^2(H,\nu)},\quad u\in D(N_2),\;v\in L^2(H,\nu),\\
\\
(u,v)_{L^2(H,\nu)}-(N^*_2v,u)_{L^2(H,\nu)},\quad u\in L^2(H,\nu),\;v\in D(N_2^*),
\end{array}\right.
$$
is a generalized Dirichlet form on
$L^2(H,\nu)$ in the sense of \cite[Definition I.4.8]{st99} with
$$
\mathcal F:=\left(D(N_2),\|N_2\cdot\|_{L^2(H,\nu)}+\|\cdot\|_{L^2(H,\nu)}\right)
$$
and with coercive part $\mathcal A$ identically equal to $0$.

We emphasize here that the theory of generalized Dirichlet forms, in
contrast to earlier versions (cf. e.g. \cite{FOT94}, \cite{MR92}),
does not require any symmetry or sectoriality of the underlying
operators. We refer to \cite{st99} for an excellent exposition. As is
well known to the experts on potential theory on $L^2$-spaces (and as is
clearly presented in \cite{st99}), the following two main ingredients
are needed.

\begin{enumerate}
\item[(a)] There exists a core $C$ of
$(N_2,D(N_2))$ which is an algebra consisting of functions having
(quasi) continuous $\nu$-versions.

\item[(b)] The capacity determined by  $(N_2,D(N_2))$ is tight.
\end{enumerate}

(a) follows from the essential $m$-dissipativity of $N_0$ on
$C_b^2(H)$ proved in the previous section, so we can take
$C:=C_b^2(H)$. This is exactly why essential $m$-dissipativity is so
important for probability theory, in particular, Markov processes.
Before we prove (b) we recall the necessary definitions.

Let
$$
G_\lambda^{(2)}:=(\lambda-N_2)^{-1},\quad \lambda>0,
$$
be the resolvent corresponding to $N_2$. A function $u\in L^2(H,\nu)$ is called
{\em $1$-excessive} if $u\ge 0$ and
$\lambda G_{1+\lambda}u\le u$ for all $\lambda>0$. For an open set $U\subset H$ define
$$
e_U:=\inf\{u\in L^2(H,\nu)|\; u \hbox{is $1$-excessive,
$u\ge 1_U$\ $\nu$-a.e.}\},
$$
(cf. \cite[Proposition III.1.7 (ii)]{st99}),
and the 1-{\em capacity} of $U$ by
$$
\mbox{\rm Cap}\;U:=\int_He_Ud\nu.
$$
(cf. \cite[Definition III.2.5 with $\varphi\equiv 1$]{st99}).
Cap is called {\em tight} if there exist increasing
compact sets $K_n,\; n\in \N$, such that for $K_n^c:=H\setminus K_n$
one has
$$
\lim_{n\to \infty}\mbox{\rm Cap}\;(K_n^c)=0.
$$
Once we have proved this, i.e., have proved (b), Theorem
 \ref{t5.1}-(i) follows from  \cite[Theorem
 IV.2.2]{st99}. Indeed, in our situation,
according to~(a) and  \cite[Proposition IV.2.1]{st99},
the requirement in \cite[Theorem IV.2.2]{st99} that
quasi-regularity holds is equivalent to
(b) and condition $D3$ in \cite[Theorem IV.2.2]{st99}.

\begin{Remark}
\label{r5.3}
{\em We mention here that in Theorem \ref{t5.1} we do not state all
facts known about $\M$; e.g., it is also proved in \cite[Theorem
IV.2.2, see also Definition IV 1.4]{st99} that all ``$\nu$-a.e.''
statements can be replaced by ``quasi everywhere'' (with respect to Cap)
statements and that

$$
x\mapsto\int_0^{+\infty}e^{-\lambda t}p_tf(x)dt
$$
is Cap-quasi-continuous. Furthermore, \cite[Theorem IV.2.2]{st99} only
claims that $\M$ has cadlag paths, but a similar proof as that in
\cite[Chapter~V, Sect.~1]{MR92} gives indeed continuous paths because
$N_2$ is a local operator.
 }
\end{Remark}

To prove (b) it is enough to find a $1$-excessive function
$u\colon\, H\to \R^+$ so that for each $n\in \N$ the level set $\{u\le n\}$
is contained in the union of a compact set $K_n\subset H$ and a $\nu$-zero set, because then $e_{K_n^c}\le
\frac1n\;u\;\nu$-a.e., hence
\begin{equation}
\label{e5.1}
\mbox{\rm Cap}\;(K_n^c)\le \frac1n\;\int_H u d\nu\to 0,\quad \mbox{\rm as}\;n\to\infty.
\end{equation}
So, the proof of Theorem \ref{t5.1}-(i) is completed by
Proposition \ref{p5.4} below, since closed balls in $L^2(D)$ are
compact in $H.$ Before we can formulate it, we need to introduce the
resolvent generated by $N_0$ on $L^1(H,\nu)$. To this end we note that
by \eqref{e2.4} $(N_0,C^2_b(H))$ is also dissipative on $L^1(H,\nu)$
(cf. e.g. \cite[Appendix~B, Lemma~1.8]{E99}), hence closable. We
recall that $(\lambda-N_0)(C^2_b(H))$ is dense in $L^2(H,\nu)$ (by the
proof of Theorem \ref{t4.1} above), hence also dense in $L^1(H,\nu)$,
so analogously $(N_1,D(N_1))$ generates a $C_0$ semigroup
$(e^{tN_1})_{t\ge 0}$ of contractions on $L^1(H,\nu)$ and we can
consider the corresponding resolvent

$$
G_\lambda^{(1)}:=(\lambda-N_1)^{-1},\quad \lambda>0,
$$
Clearly, $G_\lambda^{(1)}=G_\lambda^{(2)}$ on $(\lambda-N_0)(C^2_b(H)),$ hence
\begin{equation}
\label{e5.2}
G_\lambda^{(1)}f=G_\lambda^{(2)}f\quad\mbox{\rm for all}\;\lambda>0,\;f\in L^2(H,\nu).
\end{equation}
Define
$$
\overline{\Psi}(t):=\int_0^t\Psi(s)ds,\quad t\in \R ,
$$
and
$$
\Phi(x):=\left\{\begin{array}{l}
\int_D\overline{\Psi}(x(\xi))d\xi,\quad x\in H_ \Psi\\
\\
+\infty\quad\mbox{\rm otherwise}.
\end{array}\right.
$$
By (H3) $\overline{\Psi}$ is convex and since $r>1,$ (H3) also implies that for all $s\in \R$
\begin{align}
\label{e5.3}
0\le \frac{\kappa_0}{r(r+1)}\;|s|^{r+1}\le
\overline{\Psi}(s)
&\le  \frac{C_1}{2}\;|s|^2+
\frac{\kappa_1}{r(r+1)}\;|s|^{r+1}
\notag
\\
&\le
\Bigl[\frac{C_1}{2}\;+\Bigl(\frac{C_1}{2}\;+\frac{\kappa_1}{\kappa_0(r+1)}
\Bigr)
\;|\Psi(s)|\Bigr]\;|s|.
\end{align}
Hence, it follows by Remark \ref{r2.2} that $\Phi\in L^1(H,\nu)$.
Recall that by \eqref{e2.2}   $|\Delta\Psi|^2_H\in L^1(H,\nu).$

\begin{Proposition}
\label{p5.4}
Consider the situation of Theorem {\rm\ref{t5.1}}. Then
\begin{enumerate}
\item[{\rm(i)}]
There exists $c>0$ such that
$$
c|x|_{L^2(D)}^{r+1}\le G_1^{(1)}(\Phi+|\Delta\Psi|^2_H)(x) =:
g(x)\;(\ge 0)\;
\hbox{ for $\nu$-a.e. $x\in H.$}
$$

\item[{\rm(ii)}] The function $g^{1/2}$ is $1$-excessive.
\end{enumerate}
\end{Proposition}

For the proof of Proposition \ref{p5.4} we need the following lemma.

\begin{Lemma}
\label{l5.5}
Let $v\in C^2(H)\cap L^1(H,\nu)$ be
such that $v, |Dv|_{H^1_0}, \sup_{i\in \N}|D^2v (e_i,e_i)|$ are bounded on $H$ balls and
\begin{equation}
\label{e5.4}
\int_D\Bigl[|v(x)|\;|x|_H^2+|Dv(x)|_{H^1_0}+|Dv(x)|_{H^1_0}\;|x|_H
+ \sup_{i\in \N}|D^2v(x)(e_i,e_i)|\Bigr]\, \nu(dx)<+\infty.
\end{equation}
Then $v\in D(N_1)$ and for $\nu$-a.e. $x\in H$ one has
\begin{equation}
\label{e5.5}
N_1v(x)=\sum_{i=1}^\infty D^2v(x)(e_i,e_i)+Dv(x) (\Delta\Psi(x)).
\end{equation}
\end{Lemma}
\begin{proof}
Let $\chi\in C^\infty(\R)$ be such that $\chi'\le 0$, $0\le \chi\le 1,$ $\chi=1$ on $(-\infty,1]$
and $\chi=0$ on $(2,\infty)$. For $n\in
\N$ let
$$
\chi_n(x):=\chi\left(\frac{|x|_H^2}{n^2}\right),\quad x\in H,\quad
v_n:=\chi_n v.
$$
Then for any $x\in H$ one has
\begin{align}\label{e5.6}
Dv_n(x) &=\chi_n(x)Dv(x)+v(x)D\chi_n(x)
\notag
\\
&= 1_{\{|x|_H\le 2n\}}(x)\Bigl[\chi_n(x)Dv(x)+\frac{2}{n^2}\;v(x)
\chi'\Bigl(\frac{|x|_H^2}{n^2}\Bigr) \langle x,\cdot\rangle_H  \Bigr].
\end{align}
Likewise for $i\in \N$, $x\in H$, one has
\begin{align}\label{e5.7}
&D^2v_n(x)(e_i,e_i)
\\
&=\chi_n(x)D^2v(x)(e_i,e_i)
+v(x)D^2\chi_n(x)(e_i,e_i)+2Dv(x)(e_i)D\chi_n(x)(e_i)
\notag
\\
& =1_{\{|x|_H\le 2n\}}(x)\Bigl[\chi_n(x)D^2v(x)(e_i,e_i)
\notag
\\
&+v(x)\Bigl(\chi'\Bigl(\frac{|x|_H^2}{n^2}\Bigr)\frac{2}{n^2}\;+
\frac{4}{n^4}\;\chi''\Bigl(\frac{|x|_H^2}{n^2}\Bigr)
\langle x,e_i\rangle_H^2   \Bigr)
\notag
\\
&
+\frac{4}{n^2} \;Dv(x)(e_i)\chi'\Bigl(\frac{|x|_H^2}{n^2}\Bigr)
\langle x,e_i\rangle_H  \Bigr].\notag
\end{align}
Hence $v_n\in C^1_b(H)$. Since $|\Delta\Psi|_H\in L^2(H,\nu)$
by \eqref{e2.2} and
\begin{equation}
\label{e5.8}
\int_H|x|_H^{2r}\nu(dx)\le c_1 \int_H|x|_{L^{2r}}^{2r}\nu(dx)\le c_2 \int_H|\Delta\Psi(x)|^2_H\nu(dx)<+\infty,
\end{equation}
 (cf. Remark \ref{r2.2}), we see from \eqref{e5.6}, \eqref{e5.7} that
$v_n\to v$ and $N_0v_n$
converge to the right hand
side of \eqref{e5.5} in $L^1(H,\nu)$ as $n\to \infty$.
\end{proof}

\begin{proof}[Proof of Proposition \ref{p5.4}.]
Consider the Moreau approximation  $\Phi_\varepsilon,\;\varepsilon>0,$
of $\Phi$, i.e.,
$$
\Phi_\varepsilon(x):=\min\left\{\frac{1}{2 \varepsilon}\;\|y-x\|^2+\Phi(y)|\;y\in H   \right\},\quad x\in H.
$$
Then $\Phi_\varepsilon\in C^1(H)$, is convex and $D\Phi_\varepsilon$ is just the Yosida approximation $F_\varepsilon$
of $(\Delta\Psi,H_\Psi)$ used in \S4.
Furthermore, $\Phi_\varepsilon\uparrow \Phi$ as $\varepsilon\downarrow 0$
(cf. e.g. \cite[Proposition IV.1.8]{Sho97}).

Fix $\varepsilon,\beta>0$ and define
\begin{equation}
\label{e5.9}
\Phi_{\varepsilon,\beta}(x):=\int_H\Phi_\varepsilon(e^{\beta B}x+y)N_{\frac12\;B^{-1}(e^{2\beta B}-I)},\quad x\in H,
\end{equation}
where $B$ is as in \eqref{e4.1}. Then $\Phi_{\varepsilon,\beta}\in C^\infty(H)$, is convex and
\begin{equation}
\label{e5.10}
D_H\Phi_{\varepsilon,\beta}(x):=\Delta(D\Phi_{\varepsilon,\beta}(x))=
F_{\varepsilon,\beta}(x),\quad x\in H,
\end{equation}
with $F_{\varepsilon,\beta}$ as defined in \eqref{e4.1}. So, by the properties of $F_{\varepsilon,\beta}$ stated in \S 4 it
follows that $D^2\Phi_{\varepsilon,\beta}$ is bounded and \eqref{e4.3} implies that
\begin{equation}
\label{e5.11}
|\Phi_{\varepsilon,\beta}(x)|\le 2C_{\varepsilon}(1+|x|^2_H),\quad x\in H.
\end{equation}
By \eqref{e5.8}, \eqref{e5.11} and \eqref{e4.3} it follows that all
assumptions in Lemma \ref{l5.5} for $v:= \Phi_{\alpha,\beta}$ are
fulfilled (note that condition \eqref{e5.4} indeed holds by
\eqref{e5.8} since $r\ge 2$). Hence $\Phi_{\alpha,\beta}\in D(N_1)$
and if we denote the right hand side of \eqref{e5.5} for
$v:=\Phi_{\alpha,\beta}$ by $N_0\Phi_{\alpha,\beta}$ it follows that
for all $x\in H$
one has
\begin{equation}
\label{e5.12}
(1-N_0)\Phi_{\varepsilon,\beta}(x)\le\Phi_{\varepsilon,\beta}(x)-\langle D_H\Phi_{\varepsilon,\beta}(x),
\Delta \Psi(x) \rangle_H.
\end{equation}
Here we used that $D^2\Phi_{\varepsilon,\beta}(x)(e_i,e_i)\ge 0,\;i\in \N$, since $\Phi_{\varepsilon,\beta}$ is convex.
Since by \eqref{e4.3} one has
\begin{align*}
|\langle D_H\Phi_{\varepsilon,\beta}(x),\Delta \Psi(x)\rangle_H|
&\le C_\varepsilon(1+|x|_H)|\Delta\Psi(x)|_H
\\
&\le
C_\varepsilon(1+|x|_H)|\Psi(x)|_{H^1_0}
\end{align*}
and the right hand side is in $L^1(H,\nu)$ by \eqref{e5.8} and
\eqref{e2.2}, the right hand side of \eqref{e5.12} converges to
$\Phi_{\varepsilon}-\langle D_H\Phi_{\varepsilon,\beta}(\cdot),\Delta
\Psi(\cdot)\rangle_H$ in $L^1(H,\nu)$ as $\beta\to 0$. Applying
$G^{(1)}_1$ to \eqref{e5.12} and letting $\beta\to 0$ we then obtain
for $\nu$-a.e. $x\in H$
\begin{equation}
\label{e5.13}
 \Phi_{\varepsilon}(x)\leq G^{(1)}_1\Bigl(\Phi_{\varepsilon}(x)-\langle D_H\Phi_{\varepsilon}(x),
\Delta \Psi(x) \rangle_H\Bigr).
\end{equation}
But by \eqref{e3.6} for every $x\in H_\Psi$ one has
\begin{align*}
|\langle D_H\Phi_{\varepsilon}(x),\Delta \Psi(x)\rangle_H|
&=|\langle F_{\varepsilon}(x),F(x)\rangle_H|\\
&\le |F(x)|_{H}^2=|\Psi(x)|^2_{H^1_0}.
\end{align*}
Since $\nu(H_\Psi)=1$ and since $\Phi_{\varepsilon}+|\Psi|_{H^1_0}\in
L^1(H,\nu)$, by \eqref{e5.13} this implies that
$$
\Phi_{\varepsilon}\le G^{(1)}_1\left(\Phi_{\varepsilon}+|\Psi|^2_{H^1_0}\right)=g\quad
\nu\mbox{\rm-a.e.}
$$
Since $\Phi_{\varepsilon}\uparrow \Phi $ and $\Phi\in L^1(H,\nu)$
and  since by \eqref{e5.3} one has
$$
\Phi(x)\ge \frac{\kappa_0}{r(r+1)}\;|x|^{r+ 1}_{L^{r+1}(D)},\quad x\in
H,
$$
and $r+1\ge 2$, assertion (i) follows. To prove (ii)
fix $\lambda>0$. We note that
by the resolvent equation $\lambda G^{(1)}_{\lambda+1}g\le g,$ since $g\ge 0$. Hence
$$
\lambda G^{(1)}_{\lambda+1}g^{1/2}\le \frac{\lambda}{\lambda+1}\;
\Bigl((\lambda+1)G^{(1)}_{\lambda+1}g\Bigr)^{1/2}
=\frac{\lambda^{1/2}}{(\lambda+1)^{1/2}}\;
\Bigl(\lambda G^{(1)}_{\lambda+1}g\Bigr)^{1/2}\le g^{1/2}.
$$
So, by \eqref{e5.2} assertion (ii) follows. 
\end{proof}

The last result of this section is that in some cases the Markov processes in Theorem 5.1
can even be chosen to be strong Feller
on supp $\nu$ if $d=1$. More precisely,  consider the following condition
$$
d=1\;\;\mbox{\it and}\;C=(-\Delta)^{-\gamma}\;\;
\mbox{\it with}\;\gamma\in (1/2,1].\leqno{(C1)}
$$

\begin{Theorem}
\label{t5.5}
Assume that $(H1)$--$(H3)$ and $(C1)$
hold. Then the conservative strong Markov
process $\mathbb M$ in Theorem \ref{t5.1} can be chosen to be strong
Feller on supp~$\nu$. More precisely, its semigroup satisfies $p_tf
\in C_b(supp\;\nu)$ for all $f\in B_b(H),\;t\ge 0$, and $\lim_{t\to
0}p_tf(x)=f(x)$ for all $x\in supp\;\nu$ and all bounded Lipschitz
continuous functions $f\colon\, H\to \R$. Furthermore, $supp\;\nu$ is
an invariant set for $\mathbb M$ and Theorem {\rm\ref{t5.1}-(ii)} holds
with $\overline{H}=supp \;\nu.$

\end{Theorem}
\begin{proof}
The line of argument is exactly analogous to \cite{DPR02}. We only
mention here that the crucial estimate (4.7) in \cite{DPR02} can be
derived in the same way in our situation here. Hypotheses 1.1(i) and
1.2(i) of \cite{DPR02} are not used for this.
\end{proof}

\begin{Remark}
\begin{enumerate}[(i)]
\item
{\em
We stress that according to Theorem~\ref{th6.1} below we have that
$\mathop{supp}\nu=H$ since (C1) implies condition (H4) below.
}
\item
{\em 
For the interested reader who would like to check the details
from \cite{DPR02} for the proof of Theorem \ref{t5.5} we would like to
point out an annoying misprint in \cite[Lemma 5.6]{DPR02}. The last two
lines of its statement should be replaced by `` and for $t,\lambda>0,$
$x\mapsto \int_0^t\overline{p}_sf(x)e^{-\lambda s}ds$ is continuous on
$H_0$ ''.
}
\end{enumerate}
\end{Remark}

\section{Support of invariant measure}

In this section, we show that any measure which is the weak
limit of a sequence
of invariant probability measures $\nu_n$ corresponding to the finite
dimensional approximations has full support in the negative Sobolev
space $H:=H^{-1}(D)$ with its natural Hilbert norm $|\,\cdot\,|_{_H}$.
To this end, we obtain a uniform lower bound of $\nu_n$-measures of any
given ball in~$H$.

Let $C$ be a positive symmetric operator on $L^2(D)$.
We assume that in addition to (H1) the operator $C$ satisfies the
following condition:

\begin{itemize}
     \item[(H4)]
{\it
$\lambda_k$, $k\in\N$, in (H1) are strictly positive and there 
is a Hilbert space $E$ such that
the embedding $L^2(D)\to E$ is Hilbert--Schmidt and
$\sqrt{C}$ extends to an operator in $L(E,L^2(D))$ that will be denoted
by the same symbol.
}\end{itemize}

A typical example is $C=(-\Delta)^{-\sigma}$,
$\sigma >d/2$, and
$E=H^{-\gamma}(D)$, $\gamma >d/2$.

Let $W$ be a cylindrical Wiener process in $L^2(D)$.
Then $W$ is a continuous Wiener process with values in~$E$.
Given a
function $\Psi$ as above, we consider the
mapping $F\colon\, x\mapsto \Delta (\Psi\circ x)$ on $L^2(D)$ with values
in $H^{-2}(D)$. 

As above, let $\{e_i\}$ be the eigenbasis of the Laplacian, let
$P_n$ be the orthogonal projection in $L^2(D)$ (and also in
$H^{-1}(D)$)
 to the linear span $E_n$ of $e_1,\ldots,e_n$, and let
$F_n:=P_nF$ and $C_n:=P_n\sqrt{C}$.
We observe that
$$
\int_D \Psi\circ x(u) \Delta x(u)\, du
=-\int_D \Psi'\circ x(u) |\nabla x(u)|^2\, du \le
-\kappa \int_D |x(u)|^{r-1}| |\nabla x(u)|^2\, du
$$
for all $x\in E_n$. Therefore, on every subspace $E_n$ we have
$(F_n(x),x)_{L^2(D)}\to -\infty$ as $\|x\|_{L^2(D)}\to\infty$.
Since $F_n$ is continuous and dissipative on
$E_n$,
there is a diffusion process $\xi_n$ on $E_n$ governed
(in the strong sense) by the
stochastic differential equation
$$
d\xi_n=F_n(\xi_n)dt + C_n dW .
$$
This process has a unique invariant probability $\nu_n$.

\begin{Theorem}\label{th6.1}
Suppose that (H1)--(H4) hold and that
$1\le d\le 2(r+1)/(r-1)$.
Then any measure $\nu$ that is the limit of a weakly
convergent subsequence of $\{\nu_n\}$ has full support
in~$H$, i.e., does not vanish on nonempty open sets.
\end{Theorem}
\begin{Remark}
{\em
If $\nu:=\mu$ where $\mu$ is the solution of \eqref{e1.4} constructed in
\cite{BDPR03}, then Theorem~\ref{th6.1} applies to $\nu$.
}
\end{Remark}
\begin{proof}[Proof of Theorem~\ref{th6.1}]
Let us fix $x_0,x_1\in \bigcup_{n=1}^\infty E_n$,
$\varepsilon>0$, and
consider the deterministic equation
\begin{align}\label{6.6}
&y_n' =F_n(y_n)+C_n u_n^\varepsilon , \quad t\in [0,1],\notag
\\&
y_n(0)=x_0,
\end{align}
where $u_n^\varepsilon\in L^2(0,1;E)$ is specified below.
We consider $n\ge n_0$, where $n_0$ is such that $x_0,x_1\in E_{n_0}$.
By Lemma A.1 there is $u_n^\varepsilon\in L^2(0,1;E)$ such that
as $n\to\infty$ one has
$u_n^\varepsilon\to u^\varepsilon$ strongly in $L^2(0,1;E)$ and
\begin{equation}\label{6.8}
|y_n(1)-x_1|_{_H}\le \varepsilon .
\end{equation}
Set $D_t:= D\times (0,t)$.
Letting
$v_n^\varepsilon (t):=\int_0^t u_n^\varepsilon (s)\, ds$ we obtain
$$
\xi_n(t,x_0)-y_n(t)-\int_0^t
[F_n(\xi_n(s,x_0))-F_n(y_n(s))]\, ds =
C_nW(t)-C_n v_n^\varepsilon (t).
$$
Set $z_n(t):=\int_0^t [F_n(\xi_n(s,x_0))-F_n(y_n(s))]\, ds$. Then we
arrive at the following representation:
$$
\xi_n(t)-y_n(t)+z_n(t)=C_nW(t)-C_nv_n^\varepsilon (t).
$$
Taking the inner product in $H$ with $F_n(\xi_n(t))-F_n(y_n(t))$
and integrating in $t$ over $[0,1]$, we
obtain
\begin{multline*}
-\int_0^t \langle \xi_n(s)-y_n(s), F_n(\xi_n(s))-F_n(y_n(s))\rangle_{_H}
\, ds +\frac{1}{2}|z_n(t)|_{_H}^2
\\
=\int_0^t \langle C_nW(s)-C_nv_n^\varepsilon (s),
F_n(\xi_n(s))-F_n(y_n(s))\rangle_{_H} \, ds
\\
\le |C_nW-C_n v_n^\varepsilon|_{C([0,t]; H)}
|\Psi(\xi_n)-\Psi(y_n)|_{L^1([0,t]; H)}
\\
\le K_1 |W-v_n^\varepsilon|_{C([0,t]; E)}
|\Psi(\xi_n)-\Psi(y_n)|_{L^1([0,t]; H)} ,
\end{multline*}
where condition (H4) was employed and $K_1$ is a constant.
Generic constants will be denoted by $K$ with subindices.
Taking into account that
$$
\langle \xi_n- y_n, F_n(\xi_n)-F_n(y_n)\rangle_{_H}
=\int_D (\xi_n-y_n)(\Psi(\xi_n)-\Psi(y_n))\, du,
$$
we obtain for $t=1$
\begin{multline}\label{6.11}
-\int_{D_1} (\xi_n-y_n)(\Psi(\xi_n)-\Psi(y_n))\, du \, ds
+\frac{1}{2}|z_n(t)|_{_H}^2
\\
\le
K_1 |W-v_n^\varepsilon|_{C([0,1]; E)}
|\Psi(\xi_n)-\Psi(y_n)|_{L^1(0,1; H)} .
\end{multline}
On the other hand, by the Sobolev embedding theorem
$L^{2d/(d+2)}\subset H$ and therefore
\begin{multline}\label{new1.4}
|\Psi(\xi_n)-\Psi(y_n)|_H\le
K_2|\Psi(\xi_n)-\Psi(y_n)|_{L^{2d/(d+2)}(D)}
\\
\le
K_2\Bigl(\int_D
\Bigl[|\xi_n|^{2dr/(2+d)}+|y_n|^{2dr/(2+d)}\Bigr]\, du
\Bigr)^{(d+2)/(2d)}.
\end{multline}
Similarly to (\ref{6.11}) we have
\begin{multline*}
\int_{D_1} \xi_n \Psi(\xi_n)\, du \, ds \le
\int_{D_1} x_0\Psi(\xi_n)\, du\, ds
+K_1 |W|_{C([0,1];E)}|\Psi(\xi_n)|_{L^1(0,1;H)}
\\
\le
K_3 \int_{D_1} |x_0||\xi_n|^r\, du\, ds
+ K_4 |W|_{C([0,1];E)}
\Bigl(\int_{D_1}
|\xi_n|^{2dr/(d+2)}\, du\, ds\Bigr)^{(d+2)/(2d)}.
\end{multline*}
Since under our assumption
$2dr/(d+2)\le r+1$ we obtain
$$
\int_{D_1} |\xi_n|^{r+1}\, du\, ds \le
K_5
\Bigl(|x_0|_{L^r(D)}^r+
|W|_{C([0,1];E)}^r\Bigr).
$$
Similarly, we have by (\ref{6.6})
$$
\int_{D_1} |y_n|^{r+1}\, du\, ds \le
K_6 \Bigl(|x_0|_{L^r(D)}^r+
|v_n^\varepsilon|_{C([0,1];E)}^r\Bigr).
$$
According to (\ref{new1.4}) this yields
$$
\int_0^1 |\Psi(\xi_n)-\Psi(y_n)|_H\, ds
\le K_7
\Bigl(|x_0|_{L^r(D)}^r+
|W|_{C([0,1];E)}^r+
|v_n^\varepsilon|_{C([0,1];E)}^r\Bigr).
$$
Therefore, taking into account (\ref{6.11}) we obtain
$$
|z_n(1)|_H^2\le K_8 |W-v_n^\varepsilon|_{C([0,1];E)}
\Bigl(
|x_0|_{L^r(D)}^r+
|W|_{C([0,1];E)}^r+
|v_n^\varepsilon|_{C([0,1];E)}^r\Bigr),
$$
which along with (\ref{6.8}) gives
\begin{multline*}
|\xi_n(1,x_0)-x_1|_{_H}\le\varepsilon +
|C_nW(1)-C_nv_n^\varepsilon (1)|_{_H}
+|z_n(1)|_{_H}
\\
\le
\varepsilon +K_9 |W-v_n^\varepsilon|_{C([0,1];E)}^{1/2}
\Bigl(
|x_0|_{L^r(D)}^{r/2}+
|W-v_n^\varepsilon|_{C([0,1];E)}^{r/2}+
1\Bigr).
\end{multline*}
Therefore, for all $\alpha>0$ one has
$$
P\Bigl(|\xi_n(1,x_0)-x_1|_{_H}\ge\alpha \Bigr)
\le
P\Bigl(|W-v_n^\varepsilon|_{C([0,1];E)}^{1/2}
\bigl[|x_0|_{L^r(D)}^{r/2}+|W-v_n^\varepsilon|_{C([0,1];E)}^{r/2}
+1\bigr]\ge \gamma \Bigr),
$$
where $\gamma=(\alpha-\varepsilon)/K_9$.
Now let $\alpha=2\varepsilon$ and
let $B(x_1,\alpha)$ denote the closed
ball of radius $\alpha$ in
$H$ centered at~$x_1$. Then $B_n(x_1,\alpha)=B(x_1,\alpha)\cap E_n$
is the ball of the same radius in $E_n$  centered at $x_1$
(we recall that we deal with $n$ such that $x_1\in E_n$).
Set
$$
G_n(x_0):=
P\Bigl(|W-v_n^\varepsilon|_{C([0,1];E)}^{1/2}
\bigl[|x_0|_{L^r(D)}^{r/2}+|W-v_n^\varepsilon|_{C([0,1];E)}^{r/2}
+1\bigr]\ge \varepsilon/K_9 \Bigr).
$$
By the invariance of the measure $\nu_n$ and the previous estimate
one has
$$
\nu_n\bigl(B_n(x_1,\alpha)\bigr)=\int_{E_n}
P\Bigl(|\xi_n(1,x_0)-x_1|_{H}\le \alpha\Bigr)\, \nu_n(dx_0)
\\
\ge
\int_{E_n}
[1-G_n(x_0)]\, \nu_n(dx_0) .
 $$
Letting
$$
G(x_0):=
P\Bigl(|W-v^\varepsilon|_{C([0,1];E)}^{1/2}
\bigl[|x_0|_{L^r(D)}^{r/2}+|W-v^\varepsilon|_{C([0,1];E)}^{r/2}
+1\bigr]\ge \varepsilon/K_9 \Bigr),
$$
we have $G(x_0)=\lim\limits_{n\to\infty} G_n(x_0)$.
We recall that the measures $\nu_n$ converge weakly to $\nu$ also on
the space $L^2(D)$.
By convergence of $u_n^\varepsilon$ in $L^2(0,1;E)$ we
have $v_n^\varepsilon(t)\to \int_0^t u^\varepsilon (s)\, ds
=: v^\varepsilon$ in $C([0,1];E)$.
Therefore, the functions $G_n$ converge to $G$ uniformly
on bounded sets in $L^2(D)$. Hence
$$
\int [1-G(x_0)]\, \nu(dx_0)
=\lim\limits_{n\to\infty}\int [1-G_n(x_0)]\, \nu_n(dx_0).
$$
This yields the estimate
$$
\nu\bigl(B(x_1,\alpha)\bigr)\ge\limsup\limits_{n\to\infty}
\nu_n\bigl(B_n(x_1,\alpha)\bigr)\ge \int [1-G(x_0)]\, \nu(dx_0).
$$
It remains to observe that
$G(x_0)<1$ for every $x_0$. This follows by the fact that
$W$ is a nondegenerate Gaussian vector in $C([0,1];E)$, hence
for any $\eta>0$, one has
$P\Bigl(\sup_{t\in [0,1]}|W(t)-v^\varepsilon(t)|_E<\eta\Bigr)>0$.
\end{proof}

\vskip .3in

\centerline{{\large Appendix. Approximate controllability}}

\vskip .2in

Let $H$ be a separable Hilbert space, let $F$ be an $m$-accretive operator on $H$, and
let $B\colon\, E\to H$ be a bounded linear operator on a Hilbert space $E$ such that
$Ker (B^{*})=0$. Let $\{e_i\}$ be an orthonormal basis in $H$ and
$P_nx=\sum\limits_{i=1}^n (x,e_i)e_i$ the projection to $E_n:={\rm span}(e_1,\ldots,e_n)$.
Set $F_n:=P_nF|_{E_n}$.

Given $u\in L^2(0,1;E)$,
let us consider the following nonlinear equation:
\begin{align}\label{A1}
& y'=F y+ Bu,\quad t\in [0,T], \notag \\
&y(0)=y_0.
\end{align}
We also consider finite dimensional equations
\begin{align}\label{A3}
&y_n'=F_ny_n+ P_nBu,\quad t\in [0,T], \notag \\
&y_n(0)=P_ny_0.
\end{align}
 It was proved in \cite{BDP} that equation (\ref{A1}) is
approximately controllable, i.e., given $\varepsilon>0$ and $y_0,y_1\in \overline{D(F)}$, there
is $u\in L^2(0,1;E)$ such that
$|y(1)-y_1|_{_H}\le \varepsilon$. Here we prove a sharper result in terms of the
approximating problem (\ref{A3}).

\begin{Lemma}\label{lemA1}
Given $\varepsilon >0$ and $y_0,y_1\in \overline{D(F)}$, there exists
$u_n^\varepsilon\in L^2(0,T;E)$ such that
$|y_n(T)-P_ny_1|_{_H}\le \delta(\varepsilon)$,
$\lim\limits_{n\to\infty} u_n^\varepsilon = u^\varepsilon$ in $L^2(0,T;E)$ and
$|y^{u_\varepsilon}(T)-y_1|_{_H}\le \delta(\varepsilon)$, where
$y^u$ is the solution to {\rm(\ref{A1})}
and $\lim\limits_{\varepsilon\to 0}\delta(\varepsilon)=0$.
\end{Lemma}
\begin{proof}
It suffices to prove our claim for $y_0,y_i\in D(F)$. We fix $n$ and $\varrho>0$
and consider the differential inclusion
\begin{align}\label{A7}
& z_n'\in F_nz_n-\varrho {\rm sgn} (z_n-P_ny_1)\quad \hbox{a.e. $t\in [0,T]$}\\
& z_n(0)=P_ny_0.
\end{align}
It is known (see \cite{Ba93}) that (\ref{A7}) has a unique solution
$z_n\in W^{1,\infty}([0,T],E_n)$ and
\begin{equation}\label{A8}
z_n'(t)=F_nz_n(t)-\varrho {\rm sgn} (z_n(t)-P_ny_1)\quad
\hbox{a.e. on $[0,1]$,}
\end{equation}
where ${\rm sgn} w$ is the unit vector $w/|w|$ if $w\not=0$,
${\rm sgn} 0$ is the unit ball $\{h\in H_n\colon\, |h_n|_{_H}< 1\}$.
Therefore,
\begin{equation}\label{A9}
\frac{d}{dt}|z_n(t)-P_ny_1|_{_H}+\varrho\le
|P_nFy_1|_{_H}\quad \hbox{a.e. $t>0$.}
\end{equation}
Hence, whenever
$\varrho > |Fy_1|_{_H}+|y_0-P_ny_1|T^{-1}$, we have
$|z_n(t)-P_ny_1|=0$ for all $t\ge T$. We set
$v_n(t)\in -\varrho {\rm sgn} (z_n(t)-P_ny_1)$, where
$z_n'(t)=F_nz_n(t)+v_n(t)$ a.e. on $[0,T]$. By (\ref{A9}) we see that
$t\mapsto |z_n(t)-P_ny_1|_{_H}$ is decreasing on $[0,T]$, hence
there exists $T_0\in (0,T]$ such that
$|z_n(t)-P_ny_1|_{_H}>0$ for all $t\in [0,T_0)$ and therefore
\begin{equation}\label{A11}
v_n(t)=-\varrho \frac{z_n(t)-P_ny_1}{|z_n(t)-P_ny_1|_{_H}}\quad \hbox{for $t\in [0,T_0]$.}
\end{equation}
On the other hand, by  (\ref{A8}) we have
$$
z_n'(t)=\Bigl(F_nz_n(t)-\varrho {\rm sgn}\bigl(z_n(t)-P_ny_1\bigr)\Bigr)^{0},
$$
where $(D)^{0}$ stands for the minimal section of a set~$D$.
We have therefore
\begin{equation}\label{A12}
v_n(t)={\rm Proj}_{B(0,\varrho)} F_n(P_ny_1)\quad \hbox{for $t\in [T_0,T]$.}
\end{equation}
By (\ref{A11}) and (\ref{A12}) we conclude that as $n\to\infty$ we have
convergence $v_n\to v$ in $L^2(0,T;H)$. It is clear that
$z_n\to z$ in $C([0,T]; H)$, where
$z'(t)=Fz +v$ a.e. on $[0,T]$, $z(0)=y_0$, $z(T)=y_1$. Next, letting
$B_n:=P_nB$, we define
$u_n^\varepsilon$ to be the point where the function
$|B_n u-v_n|_{L^2(0,T;H)}^2+\varepsilon |u|_{L^2(0,T;E)}^2$ attains its minimum.
We have
\begin{equation}\label{A15}
B_n^{*}(B_nu_n^\varepsilon -v_n)+\varepsilon u_n^\varepsilon =0.
\end{equation}
Finally, we define $u^\varepsilon$ to be the point where
the function
$|Bu-v|_{L^2(0,T;H)}^2+\varepsilon |u|_{L^2(0,T;E)}^2$ attains its minimum.
We have
\begin{equation}\label{A16}
B^{*}(Bu^\varepsilon -v)+\varepsilon u^\varepsilon =0 .
\end{equation}
It follows by (\ref{A15}) and (\ref{A16}) that
$u_n^\varepsilon\to u^\varepsilon$ in $L^2(0,T;E)$ as $n\to\infty$. Moreover,
since $|Bu^\varepsilon -v|_{L^2(0,T;H)}^2+\varepsilon |u^\varepsilon|_{L^2(0,T;E)}^2\le
|v|_{L^2(0,T;H)}^2$ we have by (\ref{A16}) that
$Bu^\varepsilon -v\to 0$ weakly in $L^2(0,T;H)$ as $\varepsilon\to 0$.
Replacing $\{u^\varepsilon\}$ by a suitable sequence
of the arithmetic means of
$u^{\varepsilon_i}$ we may assume that $Bu^\varepsilon\to v$ in the norm of
$L^2(0,T;H)$. Then we see that
$$
|B_n u_n^\varepsilon -v_n|_{L^2(0,T;H)}
\le
\eta_1(1/n)
+\eta_2(\varepsilon)+C_1|u^\varepsilon -u_n^\varepsilon|_{L^2(0,T; H)}
+|B_n u^\varepsilon -B u^\varepsilon|_{L^2(0,T; H)} ,
$$
where $\eta_i(s)\to 0$ as $s\to 0$, $i=1,2$. Then we obtain
$|B_n u_n^\varepsilon -v_n|_ {L^2(0,T; H)}\le
4\eta_2(\varepsilon)=:\delta(\varepsilon)$
for all $n\ge N(\varepsilon)$.
\end{proof}

We remark that this proof remain valid if $F$ is
quasi-$m$-dissipative, i.e.,
$F+\gamma I$ is $m$-dissipative for some $\gamma>0$.
In addition, $F$ may be
multivalued.

\bigskip
 {\bf Acknowledgement}.

The second named author would like to thank the Scuola Normale
Superiore di Pisa for the hospitality and the financial support during
a very pleasant stay in Pisa when most of this work was done.
Financial support of the BiBoS-Research Centre and the
DFG-Forschergruppe "Spectral Analysis,Asymptotic Distributions,and
Stochastic Dynamics" is also gratefully acknowledged. The third named
author would like to thank the University of Bielefeld for its kind
hospitality and financial support. This work was also supported by the
research program ``Equazioni di Kolmogorov'' from the Italian
``Ministero della Ricerca Scientifica e Tecnologica''.



\begin{thebibliography}{9}


\bibitem{A86} D.G. Aronson, {\it The porous medium equation}, in
Lect. Notes Math. Vol. 1224, (A.~Fasano and al. editors), Springer,
Berlin, p.~1--46, 1986.

\bibitem{Ba93} V. Barbu, {\it Analysis and control of nonlinear
infinite dimensional systems}, Academic Press, San Diego, 1993.

\bibitem{BDP02} V. Barbu and G. Da Prato,
{\it The two phase stochastic Stefan problem},
Probab. Theory Relat. Fields, {\bf 124},  544--560, 2002.

\bibitem{BDP} V. Barbu and G. Da Prato, {\it Irreducibility of the
transition semigroup associated with the two phase stochastic Stefan
problem} (to appear).

\bibitem{BR01} V. Bogachev and M. R\"ockner, {\it Elliptic equations
for measures on infinite dimensional spaces and applications,} Probab.
Theory Relat. Fields, {\bf 120}, 445--496, 2001.

\bibitem{BDPR03} V. Bogachev, G. Da Prato and M. R\"ockner, {\it
Invariant measures of stochastic generalized porous medium equations,}
Russian Math. Dokl., {\bf 396}, 1, 2004.

\bibitem{Brezis71} H. Br\'ezis, {\it Monotonicity methods in Hilbert
spaces and some applications to nonlinear partial differential
equations}, Contributions to Nonlinear Functional Analysis,
E.~Zarantonello, ed., Academic Press, New York, 1971.

\bibitem{DP76} G. Da Prato, {\it Applications croissantes et
\'equations d'\'evolution dans les espaces de Banach,} Academic Press,
1976.

 \bibitem{DPR02} G. Da Prato and M. R\"ockner, {\it Singular
dissipative stochastic equations in Hilbert spaces,} Probab. Theory
Relat. Fields, {\bf 124}, 2, 261--303, 2002.

\bibitem{DPR03a} G. Da Prato and M. R\"ockner, {\it Invariant measures
for a stochastic porous medium equation}, BiBoS--preprint 03-07-125,
to appear in Proceedings of Conference in Honour of K. It\^o, Kyoto,
2002.

\bibitem{DPR03b} G. Da Prato and M. R\"ockner, {\it Weak solutions to
stochastic porous media equation}, BiBoS--preprint 03-07-124, to
appear in J. Evol. Equat.

\bibitem{DPZ1} G. Da Prato and J. Zabczyk, {\it Stochastic equations
in infinite dimensions}, Cambridge University Press, 1992.

\bibitem{DPZ2} G. Da Prato and J. Zabczyk, {\it Ergodicity for
infinite dimensional systems}, London Mathematical Society Lecture
Notes, {\bf 229}, Cambridge University Press, 1996.

 \bibitem{E99} A. Eberle, {\it Uniqueness and non--uniqueness of
singular diffusion operators,} Lecture Notes in Mathematics 1718, Berlin,
Springer--Verlag, 1999.

\bibitem{FOT94} M. Fukushima, Y. Oshima and M. Takeda {\it Dirichlet
forms and symmetric Markov processes,} de Gruyter, Berlin, 1994.

\bibitem{MR92} Z.M. Ma and  M. R\"ockner,  {\it Introduction
to the Theory of (Non Symmetric) Dirichlet Forms,} Springer--Verlag, 1992.

\bibitem{R99} M. R\"ockner, {\it $L^p$-analysis of finite and infinite
    dimensional diffusions}, Lecture Notes in Mathematics, {\bf 1715},
    G. Da Prato (editor), Springer-Verlag, 65-116, 1999.

\bibitem{RS03b} M. R\"ockner and Z. Sobol, {\it $L^1$--theory for the
Kolmogorov operators of stochastic generalized stochastic Burgers
equations}, Preprint 2003.

\bibitem{Sho97} R.E. Showalter, {\it Monotone operators in Banach
spaces and nonlinear partial differential equations}, Math. Surweys
and Monographs, vol. 49, AMS, Providence, 1997.

\bibitem{st99} W. Stannat, {\it The theory of generalized Dirichlet
forms and its applications in Analysis and Stochastics,} Memoirs AMS,
678, 1999.

 \bibitem{SV79} D.W. Stroock and S.R.S. Varadhan,  {\it
Multidimensional Diffusion Processes}, Springer--Verlag, 1979.
 \end{thebibliography}
 \end{document}